\documentstyle[12pt,psfig]{article}

\makeatletter
\def\@begintheorem#1#2{\it \trivlist \item[\hskip \labelsep{\bf #1\ #2.}]}
\newtheorem{teo}{Theorem}[section]

\newtheorem{rem}[teo]{Remark}
\newtheorem{lem}[teo]{Lemma}

\newtheorem{cor}[teo]{Corollary}

\newtheorem{prop}[teo]{Proposition}

\def\finedim#1{{\hfill\hbox{\enspace\fbox{\ref{#1}}}}\vspace{5pt}}
\def\dim#1{\noindent{\it Proof of} {\hspace{2pt}}\ref{#1}.}
\def\compo{\,{\scriptstyle\circ}\,}
\def\cont#1{{\rm C}^{#1}}

\font\sc=cmcsc10 scaled 1200
\def\mycaption#1{\caption{{\small #1}}}

\newfont{\Bbb}{msbm10 scaled 1200}
\def\mr{{\hbox{\Bbb R}}}

\newfont{\Got}{eufm10 scaled 1200}
\newfont{\mycal}{eusm10 scaled 1200}

\def\vecto{\!\hbox{
	\begin{picture}(10,10)
	\put(2.6,4){\line(1,0){4}}
	\put(0,0){{\hbox{\mycal X}}}
	\end{picture}}}

\def\ff{{\cal F}}

\def\gg{{\cal G}}

\def\fff{{\hbox{\mycal F}}}

\def\ristr#1{{\big|}_{#1}}
\def\hatM{\widehat{M}}
\def\d{\,{\rm d}}

\title{Branched Spines and Contact Structures on 3-manifolds}

\author{Riccardo Benedetti --- Carlo Petronio}

\begin{document}

\maketitle
 
\noindent  The study of the {\em characteristic foliation} induced on an embedded
surface is one of the main themes in 3-dimensional contact topology.  In particular
(see~\cite{giroux:conv}) one can describe exactly what foliations arise (``existence results''). Moreover
one knows that the characteristic  foliation $\ff=\fff_\xi(\Sigma)$ determines uniquely the contact
structure $\xi$ in a neighbourhood of a surface $\Sigma$ (``uniqueness results''). (By convention in this
note all surfaces are closed and oriented, all 3-manifolds are oriented, all plane fields are cooriented
and all contact structures are positive.) Moreover if $\ff$ is generic enough, so that it admits a {\em
splitting curve} $\Gamma$ (see~\cite{giroux:conv})  the geometry of $\xi$ is deeply related to the
topology of the splitting; in particular $\Gamma$ determines whether $\xi$ is {\em tight} or not near
$\Sigma$ (\cite{giroux:bourb},~\cite{giroux:criteria}).  One of the basic tools in this subject is the
so-called {\em elimination lemma}  (see~\cite{giroux:conv},~\cite{elia:palla}).

One could roughly summarize the contents of the previous paragraph as follows: {\em the study of contact
structures $\xi$ on a neighbourhood of a surface $\Sigma$ can be faithfully traslated into the study of
the ($2$-dimensional) pairs $(\Sigma,\fff_\xi(\Sigma))$}. The aim of this article is to extend this
conclusion to the case of branched surfaces. (By convention in this note  a branched surface $P$ has
singularities of generic type, i.e.  its support is a quasi-standard polyhedron, and the branching is
oriented; $S(P)$ will denote the singular set of $P$.) We have in mind in particular the significant case
where $P$ is a {\em branched standard spine} of a closed 3-manifold $\hatM$ or, more specifically, $P$ is
embedded as a {\em flow-spine} (according to Ishii's~\cite{ishii} terminology) of a flow positively
transversal to a given $\xi$ (in this case $P$ faithfully encodes the homotopy class $[\xi]$ of $\xi$ as
a plane field on $\hatM$  and will be called a {\em faithful} flow-spine for $\xi$). The theory of
branched standard spines (and in particular flow-spines) was expounded in~\cite{lnm}, and this paper is
the first development of the ideas explained there in Section 9.3. By convention, $M$ is always a compact
3-manifold bounded by $S^2$, and $\hatM$ is the corresponding closed manifold. If a spine $P$ is embedded
in $\hatM$ we always assume that the ball $B=\hatM\setminus M$ is chosen so that $M$ is a regular
neighbourhood of $P$. If $P$ is a flow-spine we also assume that the flow on $M$ positively  transversal
to $P$ extends to a ``constant'' flow on $B$  (i.e. a traversing flow which is tangent to $\partial B$
along a single {\em concave} curve, see~\cite{lnm} for more details).

Since a flow-spine which is faithful for $\xi$ already carries all the information to reconstruct
$(\hatM,[\xi])$, our initial (``uniqueness'') conjecture was that the rest of the information on the geometry of $\xi$
should be encoded by the characteristic foliation $\fff_\xi(P)$, which can be defined in a natural way.
According to Eliashberg's classification, when $\xi$ is overtwisted its geometry is determined by $[\xi]$, so
the conjecture is interesting only when $\xi$ in tight. One of the achievements of this note is to establish
the conjecture, together with a suitable ``existence'' result, under the restriction that the foliation is
``S-stable''. This means that on $S(P)$ the foliation has no singularities and the tangency points are
simple, do not lie at vertices and have index $+1$ (for a  natural definition of the index). Since we can
show that up to  $\cont{0}$-perturbation of the embedding of $P$ the foliation $\fff_\xi(P)$ is S-stable, our
result is fairly general. We also prove a weaker  ``local version'' of this result. The two most important
results of this paper are the  following: 

\vspace{.2cm}
\noindent
{\bf Theorem A.}
{\em Let $P$ be a branched standard spine of $\hatM$ and let $\ff$ be
an S-stable foliation on $P$ with isolated singularities with non-zero 
divergence. Then 
there exists a contact structure $\xi$ on $\hatM$ such that $\fff_\xi(P)=\ff$.
Given another branched standard spine $P'$ of $\hatM$ and a contact
structure $\xi'$ near $P'$ such that $(P',\fff_{\xi'}(P'))$ is
abstractly diffeomorphic to $(P,\ff)$, there exist
neighbourhoods $U$ and $U'$ of $P$ and $P'$ respectively and
a contactomorphism $\phi:(U,\xi\ristr{U})\to (U',\xi'\ristr{U'})$.}

\vspace{.2cm}

\noindent
{\bf Theorem B.}
{\em Let $(P,\ff)$ be as above and suppose the contact
structure carried by $(P,\ff)$ to be tight on a neighbourhood of $P$.
Then there exists and is unique up to isomorphism a tight contact structure $\xi$
on $\hatM$ such that $\fff_\xi(P)=\ff$. If in addition we assume that
$P$ is a flow-spine and all singularities
have positive divergence then $\xi$ can be chosen such that $P$ is faithful for $\xi$,
so $\xi$ belongs to the homotopy class carried by $P$.}
\vspace{.2cm}

These theorems summarize various more specific results, some
of which hold under the weaker assumption that $P$ is a branched surface, or
the stronger one that a certain smooth embedding of $P$ is fixed.
We warn the reader that the isomorphism $\phi$ of Theorem A does not map
$P$ to $P$ in general, therefore one cannot conclude that the ``germ''
of $\xi$ is determined by $(P,\ff)$. This is why the local result
is weaker than the tight global one.  
 
In our study of $\fff_\xi(P)$ we first establish the following facts
(we provide here informal statements which will be made precise in the 
body of the paper):
\begin{enumerate}
\item\label{ff:determines:fact} 
For any branched surface $P$ embedded in $(\hatM,\xi)$,
the germ of $\xi$ on a neighbourhood of $P$ is determined by $\fff_\xi(P)$
up to isotopies which leave $P$ invariant;
\item\label{tight:ff:determines:fact}
If $P$ is a branched standard spine embedded in $\hatM$, a {\em tight} structure $\xi$ on $\hatM$
is determined by $\fff_\xi(P)$ up to isotopies which leave $P$ invariant;  
\item\label{elim:lem:fact} 
Branched versions of the elimination lemma, in which
the separatrix is a branched leaf, hold for any branched surface $P$ in $(\hatM,\xi)$;
\item\label{normal:form:fact} For any $P$ in $(\hatM,\xi)$, the 
embedding of $P$ can be slightly $\cont0$-perturbed in such a way that 
$\fff_\xi(P)$ has a certain prescribed behaviour near vertices and is S-stable;
starting from a faithful
$P$ for $\xi$ this can be achieved with $P$ still faithful;
\item\label{all:close:fact}
For any $P$ in $(\hatM,\xi)$, if $\fff_\xi(P)$ 
is S-stable then any foliation $\cont\infty$-close to $\fff_\xi(P)$ 
is obtained by slightly $\cont\infty$-perturbing 
the embedding of $P$ in $\hatM$.
\end{enumerate}
Several of these results are proved by actually refining the arguments known in
the case of surfaces.
In particular, fact~\ref{ff:determines:fact} relies on Moser's method and a crucial
remark on the vector field generated by a branched surface according
to this method. To establish fact~\ref{tight:ff:determines:fact} we use  
Eliashberg's~\cite{elia:palla} uniqueness theorem for tight structures on the ball. 
Facts~\ref{elim:lem:fact} and~\ref{all:close:fact} extend the proofs
of Giroux~\cite{giroux:conv} using {\em ad hoc} arguments to deal with singularities,
while~\ref{normal:form:fact} uses the genuine (unbranched) elimination lemma.

The reader will note that the ``local uniqueness'' expressed by fact 1 is
stronger than stated in Theorem A. This depends on
a subtlety which is worth explaining soon,
because it shows that the results for branched surfaces are substantially
(even if not formally) different from the analogues for surfaces.
The point is that, even if the branched $\cont\infty$ structure of $P$
and the notion of $\cont\infty$ embedding in $\hatM$ are intrinsically
defined, a ``universal model relative to $P$'' of a neighbourhood of $P$ in $\hatM$,
i.e. a pair $(U,P)$ with $P\subset U$, 
does not exist, while for a surface $\Sigma$  
one can use $(\Sigma\times\mr,\Sigma\times\{0\})$. If one restricts to
standard spines, the ``absolute'' diffeomorphism type of a regular neighbourhood
of $P$ in $\hatM$ is well defined, but not the way $P$ sits in it.
This implies for instance that in~\ref{ff:determines:fact} it is not
possible to interpret $\fff_\xi(P)$ as a foliation on an 
abstract model of $P$. 

Despite what just said we can prove the weak local uniqueness 
stated in Theorem A. Namely we show that contact structures 
compatible with an abstract pair $(P,\ff)$,
with $P$ branched standard spine and S-stable $\ff$,
have isomorphic restrictions
(but the isomorphism does not map $P$ to $P$). 
A remarkable consequence of this fact is that the set of isomorphism
classes of contact structures $\xi$ on $\hatM$ such that $\fff_\xi(P)=\ff$
depends only on the abstract pair $(P,\ff)$. Another consequence is the
global uniqueness in the tight case stated in Theorem B.

Concerning constructions of contact structures our main results are:
\newcommand\genericexists{6}
\newcommand\positiveexistsclosed{7}
\newcommand\tightexistsflow{8}
\begin{enumerate}
\item[6.] If $P\subset\hatM$ is a branched standard spine,
$\ff$ is S-stable and has isolated singularities
with non-zero divergence then there exists $\xi$ on $\hatM$ 
with $\fff_\xi(P)=\ff$;
\item[7.] If $P\subset\hatM$ is a 
flow-spine, $\ff$ is S-stable and $\ff$ has isolated singularities
with positive divergence then there exists 
$\xi$ on $\hatM$ such that $\fff_\xi(P)=\ff$ and $\xi$ belongs to the homotopy
class of plane fields carried by $P$;
\item[8.] If $(P,\ff)$ is as in the previous point and carries a tight structure
on a neighbourhood of $P$
then there exists a tight structure   $\xi$ on $\hatM$ such that $P$ is faithful for 
$\xi$ and $\fff_\xi(P)=\ff$.
\end{enumerate}
Concerning fact~\genericexists, note that the condition on the divergence is necessary.
Note also that fact~\positiveexistsclosed, combined with results from~\cite{lnm}, contains a
version via branched spines of the theorem of Lutz and Martinet
according to which all homotopy classes of oriented plane fields
contain contact structures. 
Facts~\genericexists,~\positiveexistsclosed~and~\tightexistsflow~are 
proved in two steps.
The first step is to define the structure in a neighbourhood of $P$.
This is the core of our approach, and the construction
is quite subtle (much harder than the analogue
for surfaces). It relies on fact \ref{all:close:fact} stated above and
on a smooth version of the embeddability of $M$ in $P\times\mr$,
proved by Gillman and Rolfsen in a PL setting in their work on the
Zeeman conjecture~\cite{g:r:1},~\cite{g:r:2}. The second step consists in
extending the structure to the ball $\hatM\setminus M$;
in~\positiveexistsclosed~this actually requires the use of 
some of the techniques of Lutz-Martinet or of 
Eliashberg, but only on a ball, not on a general 3-manifold.
In particular one can view~\positiveexistsclosed~as a proof of the Lutz-Martinet
theorem in which, referring to the
homotopic classification of plane fields, the first (homological) 
obstruction is dealt with
by means of branched spines, and the second one (a Hopf number)
using the original approach. Fact~\tightexistsflow~can be considered as a remarkable
feature of the rigidity of tight structures.

\vspace{.2cm} \noindent {\sc Acknowledgement:} We warmly thank Emmanuel Giroux and Paolo Lisca for many
helpful and stimulating conversations. 


\section{Branched standard polyhedra,\\ embeddings and foliations}\label{smooth:def:sect}

In  this section we provide a formal definition of an 
oriented $\cont{\infty}$ branched 
surface and we discuss the notion of embedding in a 3-manifold.
Let us fix some $\cont{\infty}$ smooth function $h:\mr\to\mr$ such that
$h(x)>0$ for $x<0$ and $h(x)=0$ for $x\geq0$. In $\mr^3$ we consider the following
surfaces, all oriented so that the projection on $\mr^2\times\{0\}$ is positive:

\vspace{.2cm}

\begin{tabular}{ll}
$\Sigma_1=\mr^2\times\{0\}$ & \qquad $\Sigma_2=\{(x,y,h(x)):\ x,y\in\mr^2\}$ \\
$\Sigma_3=\{(x,y,-h(y)):\ x,y\in\mr^2\}$ & \qquad $\Sigma_4= \{ (x,y,-h(-y)):\ x,y\in\mr^2\}$ \\
\end{tabular}

\vspace{.2cm}

\noindent and define $D=\Sigma_1$, $E=\Sigma_1\cup\Sigma_2$, 
$V_+=\Sigma_1\cup\Sigma_2\cup\Sigma_3$, $V_-=\Sigma_1\cup\Sigma_2\cup\Sigma_4$. 

Let $P$ be a quasi-standard polyhedron
(a fully 2-dimensional finite polyhedron with singularities of stable nature). 
We will always assume that $P$ has a fixed ``screw-orientation'' (see~\cite{manu}:
this means that a neighbourhood of $S(P)$ can be embedded in oriented 
3-manifolds, and we
know which embeddings are positive). We will endow $E$ and $V_\pm$ with the
screw-orientation induced by $\mr^3$.

We will call 
{\em oriented branched $\cont{\infty}$ structure} on $P$ a finite collection of functions
$d_i:D\to P$, $e_j:E\to P$, $v^\pm_k:V_\pm\to P$ such that:
\begin{enumerate}
\item The union of their images covers $P$, each of them is a homeomorphism onto
an open subset of $P$, and the $r_j$'s and $v^\pm_k$'s preserve the screw-orientations;
\item Let $f,g$ be two of the maps of the family, and let $A$ be a connected
component of the domain of $f^{-1}\compo g$; let $\Sigma$ be
one of the $\Sigma_i$'s contained in the domain
of $g$, and consider the restriction of $f^{-1}\compo g$ to $A\cap\Sigma$;
then (for all choices of $f,g,A,\Sigma$) this map should take values
in one of the $\Sigma_i$'s and should be oriented and
$\cont{\infty}$ smooth as a map between surfaces.
\end{enumerate}
When $P$ is endowed with such a structure we call it a $\cont\infty$
{\em branched surface} (from now on we will always omit to specify the
{\em orientation}). If $P,\tilde P$ are branched surfaces, a map $a:P\to\tilde P$
is called a {\em diffeomorphism} if it is a homeomorphism and given any two
of the functions $f$ and $\tilde f$ which define the $\cont\infty$ structures, the
restriction of $\tilde f^{-1}\compo a\compo f$ to each of the  
$\Sigma_i$'s takes values in one of the $\Sigma_i$'s and
is smooth and oriented. 
Two structures on the same $P$ will be viewed as equivalent
if the identity is a diffeomorphism. One sees in particular that the
equivalence class of a {\em branched $\cont{\infty}$ structure} does
not depend on the particular function $h$ fixed at the beginning.

Let us recall that in~\cite{lnm} we have introduced the notion of a 
(combinatorial) branching on a quasi-standard polyhedron $P$ exactly to translate
the idea that a tangent plane should be well-defined everywhere on $P$.
We have also shown (this will be sufficient for the sequel) that 
in the oriented case a branching is just an orientation for each of
the components of $P\setminus S(P)$ such that no edge of $S(P)$ is induced the same
orientation 3 times. The following is established quite easily:

\begin{prop} 
An oriented combinatorial branching on $P$ allows to define on $P$
a structure of $\cont\infty$ branched surface, unique up to diffeomorphism.
\end{prop}

If $P$ is $\cont{\infty}$, we will call {\em smooth surface contained in $P$}
any subset locally contained in one of the $\Sigma_i$'s and open there.
A function $a:P\to\mr$ is called {\em smooth} if $a$ is smooth when restricted to
each smooth surface contained in $P$. In a similar way one defines smooth functions
from $\mr$ to $P$. We will denote in the sequel by $N$ an arbitrary (open or closed)
oriented smooth 3-manifold. An {\em embedding} of $P$
into $N$ is a map $i:P\to N$ which is a homeomorphism of $P$ into its image and is
an embedding in the usual sense when restricted to each smooth surface contained in $P$.
We will always tacitly assume that $i$ respects the screw-orientation.

Having in mind the case of genuine surfaces, an important difference arises 
when one considers embeddings of branched surfaces in 3-manifolds.
Namely, if $i_0,i_1:\Sigma\to N$ are embeddings of an oriented compact surface
into an oriented 3-manifold, then
$i_1\compo i_0^{-1}$ extends to a diffeomorphism between
tubular neighbourhoods of the images. This is false for branched
surfaces, as the following lemma already shows in dimension two.

\begin{lem}\label{not:diffeo:lemma}
Consider the maps $f,g:\mr\to\mr$ given by $f(x)=g(x)=0$ for $x\geq0$ and
$f(x)=\exp(1/x)$, $g(x)=\exp(-1/x^2)$ for $x<0$. Let $F$ (resp. $G$)
be the union of the graphs of $f$ and $-f$ (resp. $g$ and $-g$).
Then there exists no $\cont1$ diffeomorphism $\phi:\mr^2\to\mr^2$
such that $\phi(F)=G$.
\end{lem}

\dim{not:diffeo:lemma}
We only give a sketch. If $\phi$ exists then $\phi(0)=0$. Moreover $\phi$
distorts the metric in a bounded way near 0.  The two branches of
$G$ approach each other incommensurably faster than those of $F$, and
this implies that $(\partial\phi/\partial x)(x,f(x))\to+\infty$ as 
$x\to 0^-$.
\finedim{not:diffeo:lemma}

Note however that if $P$ is a standard spine then $P$ always has a
neighbourhood which is homeomorphic, and hence diffeomorphic, to
the manifold with boundary $M(P)$ defined by $P$, 
but the way $P$ sits in $M(P)$ is only
determined up to homeomorphism, not diffeomorphism (in other words,
$i_1\compo i_0^{-1}$ may not extend to a diffeomorphism). One
way to overcome this difficulty, which we will sometimes refer to in this paper, is
to go back to the beginning of the section, choose a definite function
$h$ and require that the embedding of $P$ into $N$ should extend, in
each of the local models $D,E,V_\pm$, to an embedding of $\mr^3$.
This will be called an {\em $h$-embedding} of $P$ in $N$.
It is easily checked that indeed if $i_0,i_1$ are $h$-embeddings
for the same $h$ then $i_1\compo i_0^{-1}$ extends to a diffeomorphism
between neighbourhoods.

We will now introduce  another notion of embedding, which builds
on the results of~\cite{lnm} and relates branched polyhedra to contact structures. 
We consider a closed manifold $\hatM=M\cup B$. Let us recall that certain 
standard spines $P$ (called flow-spines)
of $M$ can be endowed by an oriented branching such that
there exists $v\in\vecto(\hatM)$ with $v$ 
positively transversal to $P$ and $(B,v)$ diffeomorphic to a 
constant field on $B^3$. Moreover $P$ determines the homotopy
class of $v$, and all classes are carried by some $P$.
Now consider a contact structure $\xi$ on $\hatM$. We will say that
$P$ is embedded as a {\em faithful flow-spine for $\xi$} if
there exists $v$ as above which moreover is positively
transversal to $\xi$. Note that in this case $P$ encodes
the homotopy class of $\xi$.

We remark now that several differentiable notions,
like tangent vectors, foliations and differential forms can be defined
in an obvious way for branched surfaces. Moreover, if $P$ 
is a branched surface embedded in a contact manifold $(N,\xi)$ then
the {\em characteristic foliation} $\fff_\xi(P)$ induced by $\xi$ on $P$
is well-defined. The following two lemmas are easy.

\begin{lem}
Let $P$ be embedded in $(N,\xi)$. Then up to a $\cont\infty$ small
perturbation of the embedding we can assume that $\fff_\xi(P)$ has isolated
singularities away from $S(P)$ and isolated simple tangency points to
$S(P)$ away from vertices. Starting from an $h$-embedding for some $h$,
or from an embedding of a faithful flow-spine, or both, one can find 
the new embedding with the same properties.
\end{lem}

\begin{lem}
\begin{enumerate}
\item If $P$ is a branched surface in $(N,\xi)$ and $\fff_\xi(P)$
is as in the previous lemma then the singularities of $\ff(P)$ have non-zero divergence.
\item If $P$ is a faithful flow-spine for $\xi$ on $\hatM$ and $\fff_\xi(P)$
is as in the previous lemma then the singularities of $\fff_\xi(P)$ have positive
divergence.
\end{enumerate}
\end{lem}

Concerning the second fact, we will always use in this note the
conventions of Giroux~\cite{giroux:bourb} on orientations. In particular, a
singularity $p$ of $\fff_\xi(P)$ has positive divergence if and only if the
orientations of $\xi$ and $P$ coincide at $p$.

\begin{figure}
    \centerline{\psfig{file=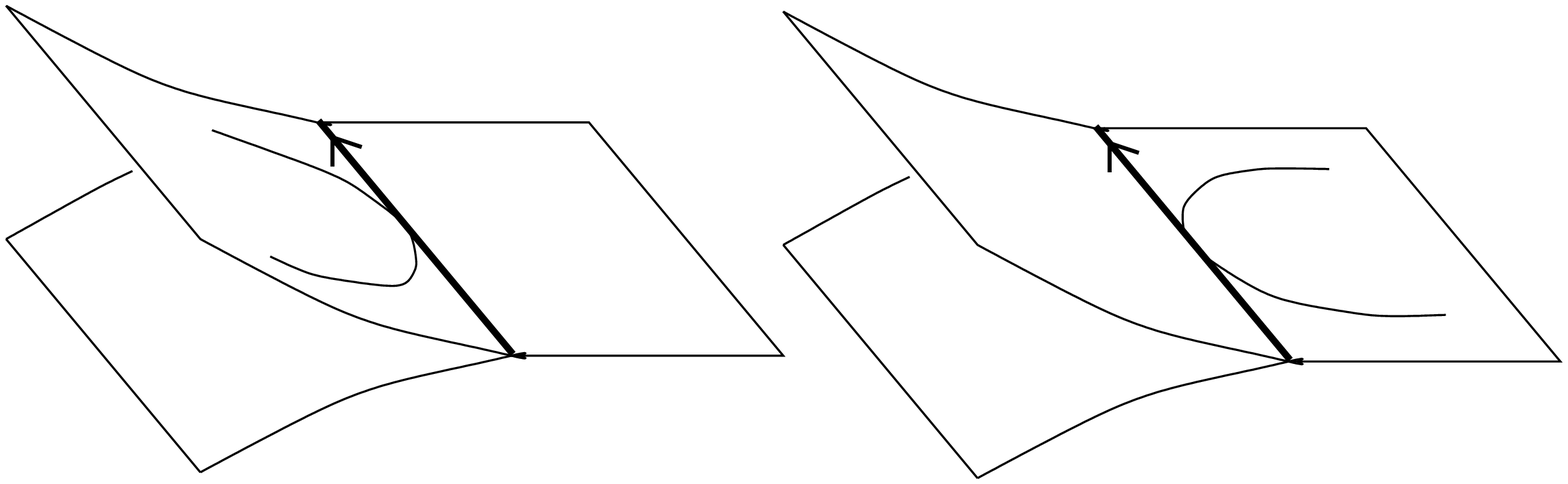,width=7truecm}}
    \mycaption{\label{def:index}
    Simple tangecies of index $+1$ and $-1$ respectively}
    \end{figure}

\vspace{.2cm}
We conclude this section with a result which justifies the definition (given
in the introduction) of S-stable foliation on a branched surface. We first remark that 
if $p$ is a simple tangency point between $S(P)$ and an oriented foliation $\ff$ on $P$,
the index of $p$ is naturally defined, as shown in Fig.~\ref{def:index}.
Recall that $\ff$ on $P$ is called S-stable if it is non-singular along $S(P)$
and has simple tangencies of index $+1$ to $S(P)$ away from $V(P)$. The following
result implies in particular that the germ of such an $\ff$ along $S(P)$ is stable
under $\cont{\infty}$-perturbation, whence the name.

\begin{prop}\label{S-stable:prop}
Let $\ff$ be an S-stable foliation on $P$. If $\ff'$ is sufficiently $\cont{\infty}$-close to
$\ff$ then there exists a regular neighbourhood $U$ of $S(P)$ and a diffeomorphism
$\varphi:P\to P$ such that $\varphi(U)=U$, $\varphi_*(\ff\ristr{U})=\ff'\ristr{U}$ and
$\varphi_*(\ff)$ is everywhere $\cont{\infty}$-close to $\ff'$.
\end{prop}

\dim{S-stable:prop}
We establish some preliminary facts on planar foliations. All foliations and curves
are $\cont{\infty}$-smooth.

{\em Claim 1. Let $\ff$ be a non-singular foliation near $0\in\mr^2$ and let $\gamma$
be a curve with a simple tangency to $\ff$ at $0$. Up to diffeomorphism we can assume
that, near $0$, $\ff$ is horizontal and $\gamma$ is the curve $t\mapsto(t,t^2)$.}
Of course we can assume $\ff$ to be horizontal. Then $\gamma$ is the graph of a function 
$f:(-\varepsilon,\varepsilon)\to\mr$ with $f(0)=f'(0)=0$ and $f''(0)>0$. So $f(x)=x^2\cdot g(x)$
with $g(0)>0$. Therefore $k=\sqrt{g}$ is smooth, and the required diffeomorphism is
$(x,y)\mapsto (x\cdot k(x),y)$.

\begin{figure}
    \centerline{\psfig{file=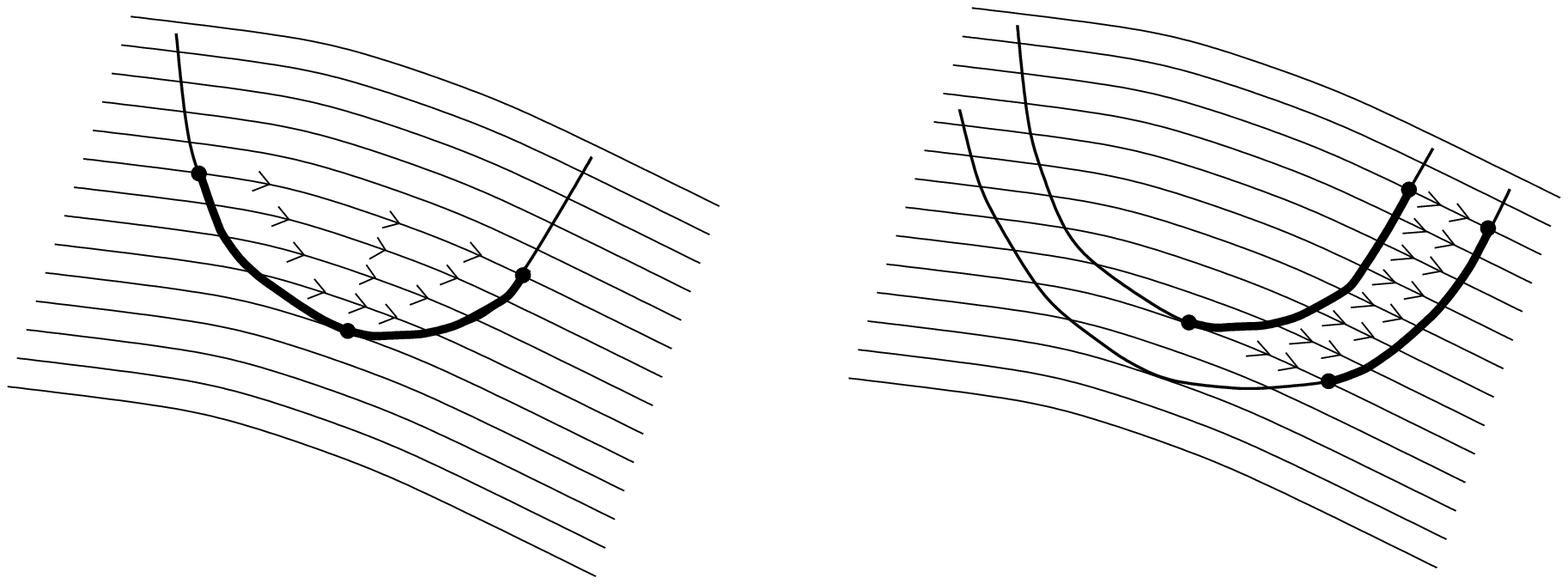,width=10truecm}}
    \mycaption{\label{alongleaf}
    Functions defined following the leaves}
    \end{figure}

{\em Claim 2. Let $(\ff,\gamma)$ be a pair as in claim 1, and consider the  
involution $g:\tilde\gamma\to\tilde\gamma$ on a closed subinterval $\tilde\gamma$ of $\gamma$ defined
as suggested on the left of Fig.~\ref{alongleaf}. Then $g$ is a diffeomorphism.}
This is obvious because in the universal model given by 
claim 1 the involution is just $t\mapsto -t$.
    
{\em Claim 3. Let $(\ff_i,\gamma_i)$, $i=0,1$, be pairs as in claim 1,
and let $\delta_i$ be a curve $\cont{\infty}$-close to $\gamma_i$
which meets twice the leaf of $\ff_i$ through 0. Consider the
map $f_i:\tilde\gamma_i\to\tilde\delta_i$ between closed subintervals of
$\gamma_i$ and $\delta_i$, defined as  
suggested on the right of Fig.~\ref{alongleaf}.
Then $f_1^{-1}\compo f_0$ is a diffeomorphism (whereas $f_0$ and $f_1$ are not).}
Claim 1 implies that we can assume that $(\ff_0,\gamma_0)=(\ff_1,\gamma_1)$,
and the conclusion easily follows.

{\em Conclusion.} We will denote by $p_i$ (resp. $p'_i$) the tangency points of $\ff$ 
(resp. $\ff'$) to $S(P)$. Note that $p'_i$ also has index +1 and is close to $p_i$.
We choose a regular neighbourhood $U$ of $S(P)$ whose boundary
is very close and almost parallel to $S(P)$, except near
vertices where it turns smoothly. For both $\ff$ and $\ff'$ the tangency points 
to $\partial U$ have the following qualitative description: there are exactly two
for each vertex of $P$ and exactly one for each $p_i$ (or $p'_i$).
For each $p_i$ we select a neighbourhood $A_i$
which is bounded by two Y-shaped leaves of $\ff$ and three segments which lie in
$\partial U$. We do the same for $p'_i$, taking $A'_i$ to be almost
identical to $A_i$, with $p_i\in A'_i$ and
$p'_i\in A_i$.

It is now quite easy to construct a diffeomorphism 
$\varphi:S(P)\setminus\bigcup_iA'_i\to S(P)\setminus\bigcup_iA_i$
such that following the leaves of $\ff'$ and $\ff$ we get
a diffeomorphism $\varphi:U\setminus\bigcup_iA'_i\to U\setminus\bigcup_iA_i$
which transforms $\ff'$ to $\ff$ (a little care has to be taken for the choice of
$\varphi$ near vertices). Closeness of diffeomorphisms to the identity will be
easy in each step of the present proof, and will not be explicitly mentioned.

\begin{figure}
    \centerline{\psfig{file=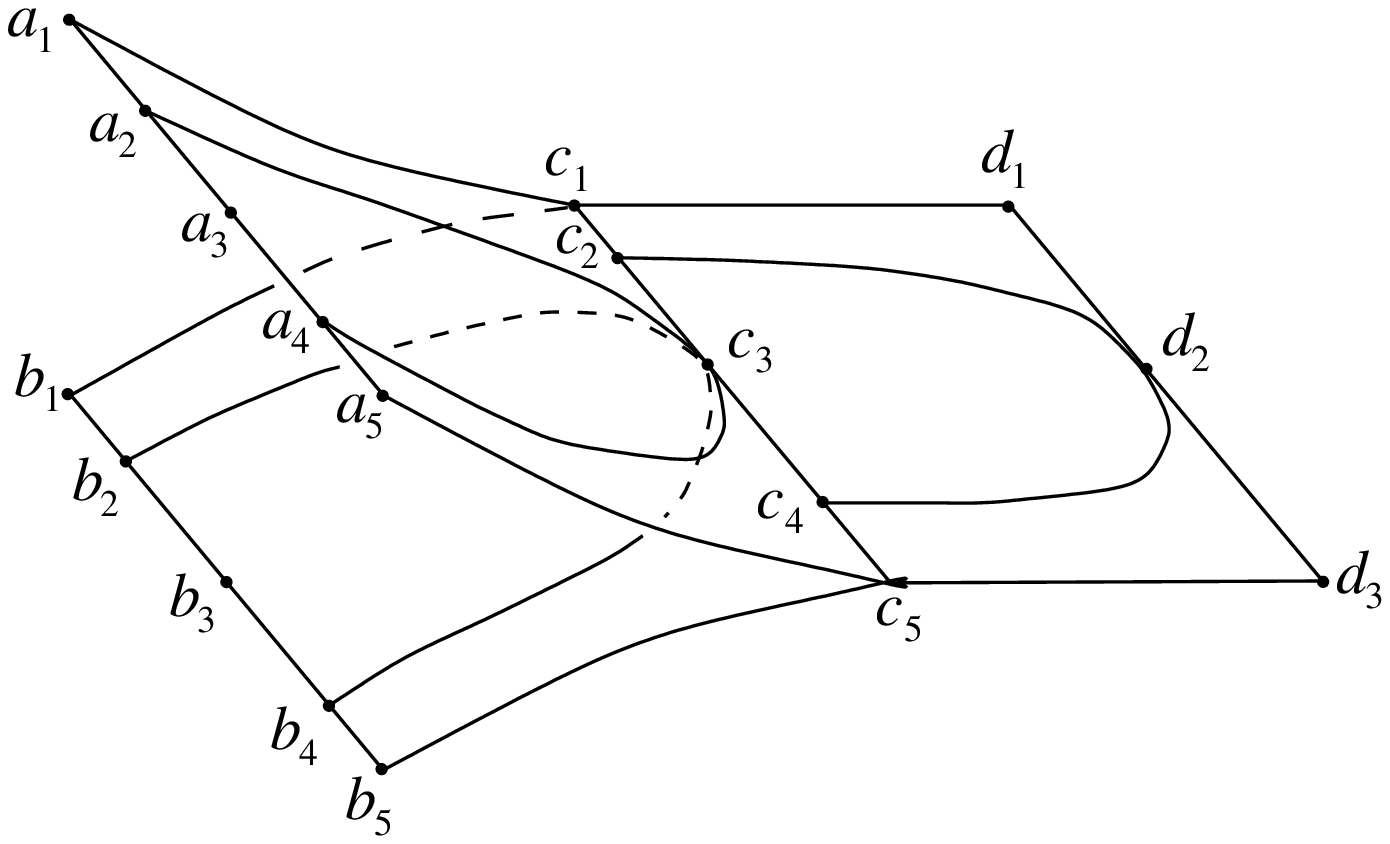,width=8.5truecm}}
    \mycaption{\label{goodtang}
    Local aspect of a singularity of index $+1$}
    \end{figure}

We will now extend $\varphi$ mapping $A'_i$ to $A_i$.
To do this we note that $A_i$ and $A'_i$ both have a 
description as shown in Fig.~\ref{goodtang}.
Following the leaves of $\ff$ we obtain various homeomorphisms between closed intervals,
for instance $[a_1,a_2]\to[c_1,c_3]$, $[c_1,c_2]\to[d_1,d_2]$ \dots, and some diffeomorphisms,
for instance $[a_2,a_3]\to[a_4,a_3]$, $[c_2,c_3]\to[c_4,c_3]$ \dots . We do the same for $\ff'$.
Now to extend $\varphi$ to a map $A'_i\to A_i$ we first extend it to $[d'_1,d'_3]\to[d_1,d_3]$
with $d'_2$ mapped to $d_2$. Following the leaves we obtain the extension of $\varphi$
to the planar quadrilaterals $(d'_1,d'_2,c'_2,c'_1)\to(d_1,d_2,c_2,c_1)$. The next step
is to further extend to $[c'_2,c'_3]\to[c_2,c_3]$ and again follow the leaves to extend to
$(c'_2,c'_3,c'_4,d'_2) \to (c_2,c_3,c_4,d_2)$. We proceed in a similar way to construct
the diffeomorphism $A'_i\to A_i$, repeatedly using the claims 2 and 3 above to show
that indeed we have a diffeomorphism and not only a homeomorphism.

The last step consists in extending $\varphi$ to the whole of $P$. Since we do not
require the foliation to be preserved, this can be done in any arbitrary way.
\finedim{S-stable:prop}

\begin{rem}\label{non-stable}
{\em The previous result does not hold if $\ff$ has a simple tangency $p$ to $S(P)$ of index $-1$.
This is because following the leaves on the two branches on the left of $S(P)$ we 
construct two germs at $p$ of involutions defined on the singular edge, and
the composition of these involutions is, up to conjugation, a non-constant invariant.}
\end{rem}

\section{Uniqueness results}\label{uniqueness:sect}

In this section we establish facts~\ref{ff:determines:fact} 
and~\ref{tight:ff:determines:fact} from the introduction. Note that the latter implies
the uniqueness part in Theorem B.

\begin{teo}\label{foliation:determines}
Let $P$ be a branched surface embedded
in $N$ and let $\xi_0,\xi_1$ be contact structures on $N$
such that $\fff_{\xi_0}(P)=\fff_{\xi_1}(P)$. Then there exist a neighbourhood $U_0$ of
$P$ and a smooth function $\phi:U_0\times[0,1]\to N$ such that:
\begin{enumerate}
\item For all $t$ the map $\phi_t=\phi(\,\cdot\,,t)$ is a diffeomorphism
of $U_0$ onto an open neighbourhood $U_t$ of $P$; moreover
$\phi_t(P)=P$ and $\phi_0={\rm id}$;
\item $\xi_0\ristr{U_0}=\phi_1^*(\xi_1\ristr{U_1})$ and 
$\fff_{(\phi_t)_*(\xi_0|_{U_0})}(P)=\fff_{\xi_0}(P)=\fff_{\xi_1}(P)$ for all $t$.
\end{enumerate}
\end{teo}

\dim{foliation:determines}
By simplicity of notation we will assume that the regions of $P$
(the components of $P\setminus S(P)$) have closure homeomorphic
to the closed disc. In general, when the 
regions are not discs or there are self-adjacencies, we would have
to slightly modify the proof by cutting the regions into portions.
Now we can parametrize the regions of $P$ by maps $f_i:Q_i\to P$ where:
\begin{enumerate}
\item $Q_i$ is either a closed disc or a regular $n$-gon for some $n$
(for $n=2$ we define the bigon as $\{(x,y):\ x^2+(y-1)^2\leq2,\ x^2+(y+1)^2\leq2\}$);
\item $f_i$ is a homeomorphism onto a closed region of $P$ and it extends to a 
smooth embedding $\tilde f_i$ of $\mr^2$ into $N$.
\end{enumerate}

Now let us consider on $N$ a Riemannian metric, and let us note that for all
$p\in P$ the positive unit normal  $\nu(p)$ to $P$ in $p$ is well-defined.
Moreover $\nu\compo f_i$ extends to a smooth function defined on $\mr^2$.
Now, up to modifying the $\tilde f_i$'s without changing the $f_i$'s,
we can find $\varepsilon>0$ such that the map
$$\tilde F_i:\mr^2\times(-\varepsilon,\varepsilon)\ni(x,y,t)\mapsto 
\tilde f_i (x,y)+t\cdot \nu(\tilde f_i(x,y))$$
is an embedding. Up to a change of scale we assume that $\varepsilon=\infty$.
Let us define $F_i$ as the restriction of $\tilde F_i$ to $Q_i\times\mr$.
Even if the $F_i$'s are defined on subsets of $\mr^3$ which are not open,
we can view them as charts, because they extend to diffeomorphisms.
Note that by the very construction each ``coordinate change'' $F_i^{-1}\compo F_j$
is the identity on the last coordinate. 
The condition that $P$ is smoothly embedded implies that 
the images of the $F_i$'s cover a neighbourhood of $P$.

The rest of the proof follows quite closely the argument yielding the same
result for surfaces, so we omit computations and confine ourselves to a
description of the various steps.
Let $F_i^*(\xi_j)$ be defined by a form 
$\alpha^{(i)}_j=\beta^{(i)}_j+u^{(i)}_j\d t$. Since $F_i^*(\xi_0)$ and
$F_i^*(\xi_1)$ induce the same characteristic foliation on $Q\times\{0\}$,
up to multiplying one of the equations by a scalar function we can assume that
$\beta^{(i)}_0=\beta^{(i)}_1$ for $t=0$. Note that the scalar function is
determined on the image of the various $F_i$'s, but one easily sees
that it glues up to a smooth function on a neighbourhood of $P$.

Now we define $\xi_s=(1-s)\xi_0+s\xi_1$. If one considers the form 
$\alpha_s=F_i^*(\xi_s)=(1-s)\alpha_0+s\alpha_1$ and computes 
$\alpha_s\wedge\d\alpha_s$, using the condition $\beta^{(i)}_0=\beta^{(i)}_1$ one sees
that this 3-form is positive for $t=0$, and hence for small enough $t$.
Since there are finitely many charts $F_i$ one deduces that $\{\xi_s\}$
is a homotopy of contact structures on some neighbourhood of $P$.

\begin{figure}
    \centerline{\psfig{file=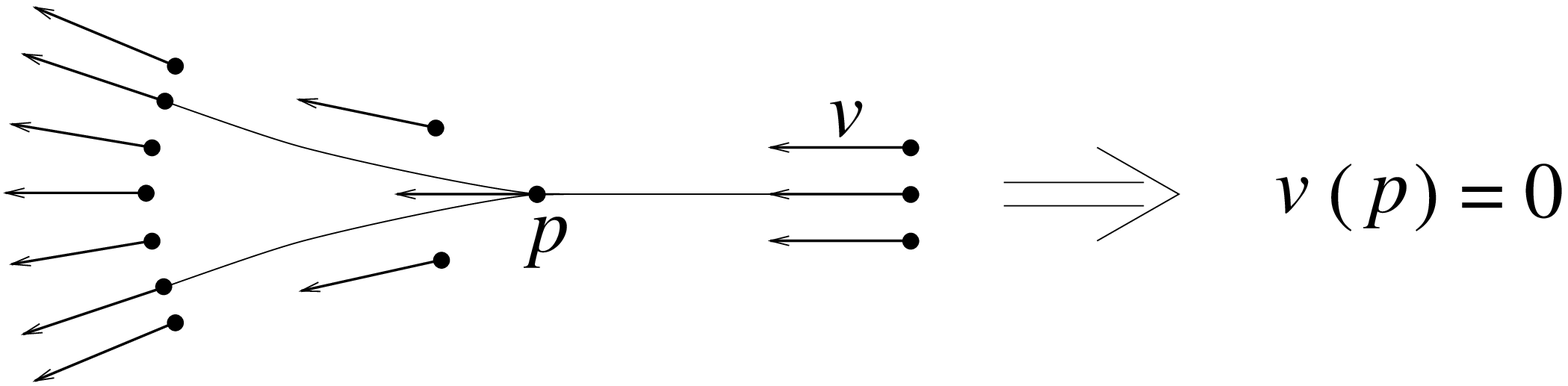,width=10truecm}}
    \mycaption{\label{viszero}
    A vanishing result for vector fields}
    \end{figure}

The next step consists in applying Moser's method. We know that 
a contact homotopy $\xi_s$ yields a time-depending vector field
$v_s$ integrating which (when possible) one conjugates $\xi_0$ to $\xi_1$.
If ones looks closely at the definition of $v_s$, one sees that if
$\Sigma$ is a surface and all the $\xi_s$'s induce on $\Sigma$ the same
characteristic foliation, then $v_s$ is always parallel to the vector 
field which directs this characteristic foliation. In our setting, since $P$
can be seen as a union of surfaces, this implies that $v_s$ is always tangent
to $P$. As a consequence, we deduce that {\em along $S(P)$ the vector
field $v_s$ is null or tangent to $S(P)$} (so in particular it is null
at vertices). One can easily verify this fact directly,
using some model of the
embedding near $S(P)$. An indirect proof can be 
obtained by contradiction as suggested by Fig.~\ref{viszero}.
If $v$ is at least $\cont1$ then it defines a dynamical system
which, as the figure shows, is non-deterministic: a contradiction.

It follows from above that the vector field $v_s$ can be integrated up
to time 1 on $P$ and the resulting map leaves $P$ (more precisely, all
its regions) invariant. Pretty much as in the case of surfaces, this
implies that $v_s$ can be integrated up to time 1 also on some neighbourhood
of $P$, yielding a diffeomorphism on some other neighbourhood. For a
formal proof, one considers the expressions 
$v^{(i)}_s=w_s^{(i)}+r_s^{(i)}\cdot\partial/\partial t$, chooses
constants $\varepsilon,\delta$ such that $|r_s^{(i)}(x,t)|\leq\delta\cdot t$
for $|t|\leq\varepsilon$, and shows by an a priori estimate on the
solution of a Cauchy problem that starting (in some chart) with a
$t$-coordinate less than $\varepsilon\cdot\exp(-\delta)$,
time 1 is reached and (in some other chart) the $t$-coordinate
is less than $\delta$. This implies the conclusion.
\finedim{foliation:determines}

\begin{prop}\label{tight:determines}
Let $P$ be a branched spine embedded in $\hatM$ and let
$\xi_0,\xi_1$ be tight structures on $\hatM$ such that
$\fff_{\xi_0}(P)=\fff_{\xi_1}(P)$. Then there exists an isotopy
between $\xi_0$ and $\xi_1$ which 
leaves $P$ invariant and preserves $\ff$.
\end{prop}

\dim{tight:determines}
Using the notations of the statement of Theorem~\ref{foliation:determines}
we can assume that $U_t$ is a regular neighbourhood of $P$ whose
boundary $S_t$ is an embedded 2-sphere. We can also assume that $\xi_t$
is defined on a neighbourhood $V_t$ of the closure of
$U_t$, and that the foliation induced by $\xi_t$ on $S_t$ is the trivial
one, with one source and one sink and no saddles or cycles.
The last condition can be imposed to $S_0$ according to~\cite{giroux:bourb},
and is then automatic for all $S_t$'s. Now we can extend each $\xi_t$ to
$\hatM$ by identifying the ball $\hatM\setminus U_t$ with the unit ball
in the standard tight structure on $\mr^3$. Eliashberg's uniqueness 
theorem~\cite{elia:palla} for tight structures on the ball implies that
the resulting family $\xi_t$ can be assumed to be 
continuous. In other words, $\xi_t$ is a contact homotopy. The conclusion
now follows using Gray's theorem (and its proof: we need to
note that $(P,\ff)$ is invariant under the flow generated by
Moser's method).
\finedim{tight:determines}

\begin{rem}
{\em The above result holds under the a priori weaker assumption
that $\xi_0,\xi_1$ should be tight on $\hatM\setminus P$.}
\end{rem}

\section{Modifying the characteristic foliation}\label{modifying:sect}

In this section we provide
formal statements and proofs of facts~\ref{elim:lem:fact},~\ref{normal:form:fact}
and~\ref{all:close:fact}. We start 
with an easy fact which answers a natural question and will be useful later.

\begin{prop}\label{homotopic:faithful}
Let $\xi$ be a plane field on $\hatM$ and let $P\subset\hatM$ be a flow-spine which carries
the homotopy class of the vector field positively transversal to $\xi$. Then $P$
can be isotoped to be faithful for $\xi$.
\end{prop}

\dim{homotopic:faithful}
Let $v$ be a vector field positively transversal to $\xi$. Then by~\cite{lnm}
(or~\cite{ishii}) there exists a flow-spine $Q$ for $v$, i.e. one which is faithful
for $\xi$. Now by~\cite{lnm} $Q$ and $P$ are related by a sequence 
standard sliding moves. If we realize the sequence of moves within $M$ we can
require that $v$ always remains positively transversal. The result is a spine $P'$
isomorphic to $P$ and faithful for $\xi$. Since the complements of $P$ and $P'$
are balls, $P$ and $P'$ are isotopic.
\finedim{homotopic:faithful}
 
We proceed now with a branched version of the elimination lemma.

\begin{teo}\label{branched:elim:state}
Let $P$ be a branched surface embedded in a contact $3$-manifold $(N,\xi)$.
Then the qualitative local modifications of $\fff_\xi(P)$ shown in
Fig.~\ref{elimlem1} and Fig.~\ref{elimlem2}~can 
    \begin{figure}
    \centerline{\psfig{file=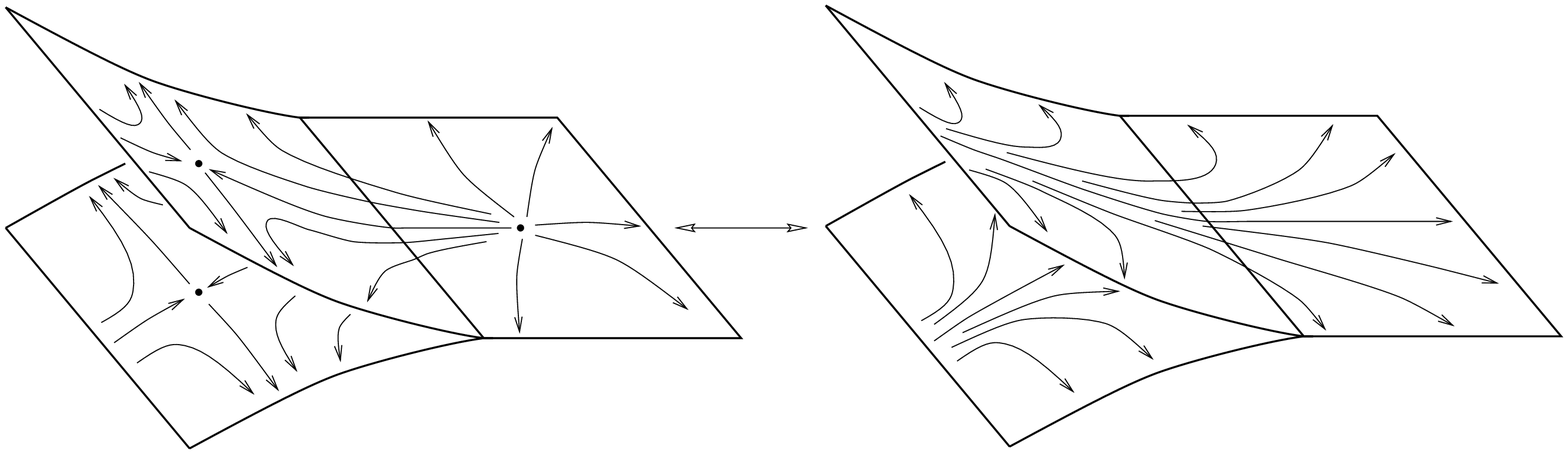,width=14truecm}}
    \mycaption{\label{elimlem1}
    Branched elimination lemma}
    \end{figure}
    \begin{figure}
    \centerline{\psfig{file=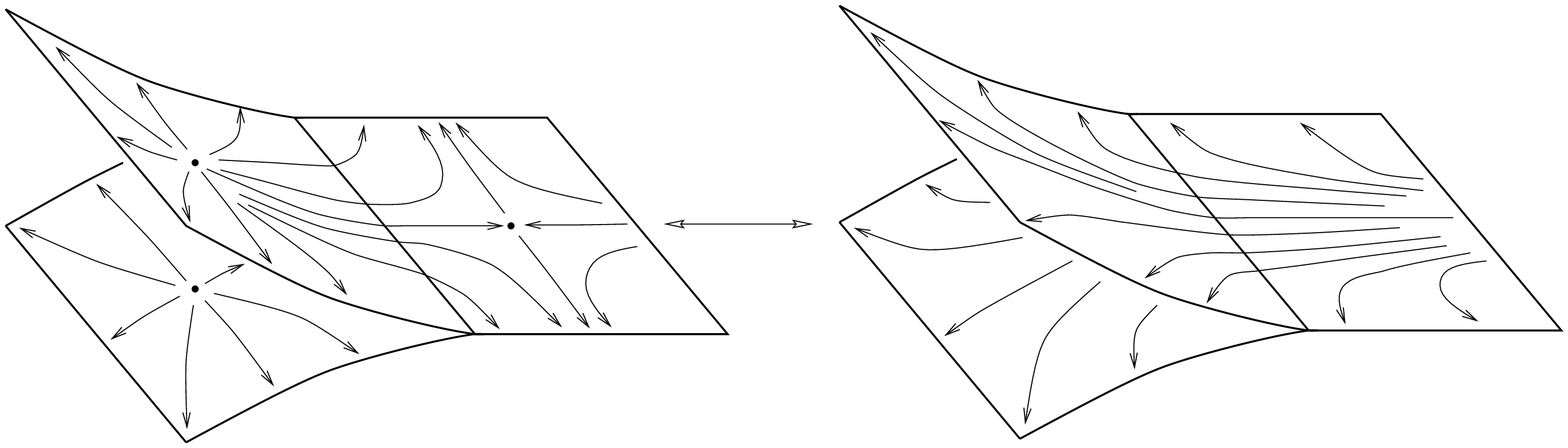,width=14truecm}}
    \mycaption{\label{elimlem2}
    Branched elimination lemma}
    \end{figure}
be achieved by $\cont0$-small perturbations
of the embedding, provided in both figures 
the saddle points are assumed to be positive. 
If we start with an $h$-embedding for some $h$ (or a faithful embedding 
of a flow-spine for $\xi$, or both) then the perturbed embedding
has the same properties. A similar result, except for the
case of faithfully embedded flow-spines, holds for negative 
rather than positive singularities.
\end{teo}

\begin{rem} {\em Just as in the case of ordinary surfaces, the statement of the elimination
lemma must be understood with some care. Namely, if we denote by $\ff_0$ and $\ff_1$
the foliations shown in  Fig.~\ref{elimlem1} (or~\ref{elimlem2}) on the same abstract branched
rectangle $R$, the result should be formally stated as follows. {\em If $P$ is an abstract branched
surface, $\{i_0,i_1\}=\{0,1\}$, $j_0:P\hookrightarrow N$ and $\varphi_0:R\hookrightarrow P$
satisfy $\fff_\xi(j_0(\varphi_0(R)))=j_0(\varphi_0(\ff_{i_0}))$, then {\em there exist}
$\varphi_1:R\hookrightarrow P$ with $\varphi_1(R)=\varphi_0(R)$ and $j_1:P\hookrightarrow N$
with $j_1=j_0$ outside $\varphi_0(R)$ and $\fff_\xi(j_1(\varphi_1(R)))=j_1(\varphi_1(\ff_{i_1}))$.}
The point we are making here is that there are various ways to fit $\ff_{i_1}$ instead of $\ff_{i_0}$
on $\varphi_0(R)$; these various ways give inequivalent global foliations on $P$, and we cannot prescribe
a priori which one will arise. A na\"\i ve statement of the elimination lemma, in which one
replaces ``$\exists\varphi_1$'' by ``$\forall\varphi_1$'', is easily seen to lead to
contradictions already in the unbranched case. For instance one could construct an overtwisted
disc around any focus, contradicting existence of tight structures.}
\end{rem}

\dim{branched:elim:state}
We will assume that the reader is familiar with the details of the 
proof of this result
for surfaces, as exposed for instance in~\cite{audin}. We will first
refer to the case of embeddings without further properties.

\begin{figure}
    \centerline{\psfig{file=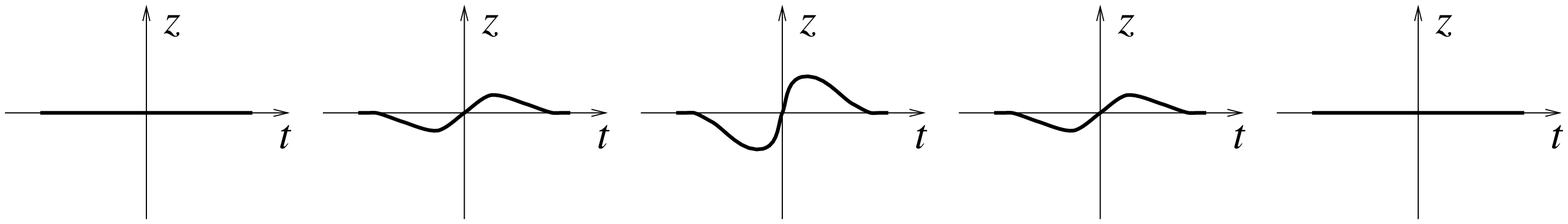,width=15.5truecm}}
    \mycaption{\label{zetafunz}
    The function $\zeta$}
    \end{figure}

We can imagine the portion of $P$ shown in the figures as the
union of two smooth rectangles $R_+$ and $R_-$ which share a square. 
The proof of the elimination lemma (in both directions) 
for each of $R_\pm$ viewed by itself would go as
follows. We first parametrize a neighbourhood of $R_\pm$ as $\mr^3_{(u,t,z)}$
so that $R_\pm$ is the $(u,t)$-plane, the singularities and 
both separatrices of the saddle lie on the
the $u$-axis and the expression of $\xi$
in coordinates $(u,t,z)$ satisfies certain properties. 
Now the new $R_\pm$ is the graph of a function
$\zeta_\pm:\mr^2_{(u,t)}\to\mr_z$ qualitatively described in Fig.~\ref{zetafunz},
where we show cross-sections $\{t=t_0\}$  and $t_0$ increases from left to right.

The idea is now to find compatible parametrizations for
$R_+$ and $R_-$, so that the union of the two perturbed rectangles
gives a perturbed copy of $P$. This can be done quite easily, one
only needs to be careful and take the common square to be the
half-plane $u\geq0$ in both parametrizations.

If we start with an $h$-embedding for some $h$ then of course we can
assume that we end up with an $h$-embedding: the qualitative
picture of the foliation is insensitive to a $\cont\infty$-small modification
of the embedding near $S(P)$.

Let us now prove that if we start with a faithful flow-spine for $\xi$ then
we end up with a faithful flow-spine. So, let us denote by 
$P_0$ and $P_1$ the initial and final embeddings of $P$, and let us
assume that there exists a vector field $v_0$ positively transversal to
both $\xi$ and $P_0$. We want to modify $v_0$ locally to a flow $v_1$
positively transversal to $\xi$ and $P_1$. We will actually do this
separately for each of the rectangles $R_\pm^{(0)}$ and $R_\pm^{(1)}$
which cover $P_0$ and $P_1$ respectively: as above, the choice of
compatible coordinates gives the desired result.

So, we refer to coordinates $(u,t,z)$ as above. Let us first remark that 
$\partial/\partial u$ is tangent to both $R_\pm$ and 
$\xi$ at the points $(u,0,0)$. So we can get rid of the $u$-coordinate of
$v_0$, and actually assume (at least in the zone affected by the 
modification) that $v_0$ is independent of $t$ and $z$.
Therefore we can concentrate on one of the planes $\{u=u_0\}$
and show how to construct $v_1$ there. Up to rotations and dilations, 
since $v_0$ is constant, we can assume that $v_0=\partial/\partial z$ 
on the plane. Now we note that $\fff_\xi(\{u=u_0\})$ is non-singular, and,
since $v_0$ is transversal to $\xi$, this foliation is transversal
to the vertical lines parallel to the $z$-axis. Therefore, up
to a change of chart of the form $(t,z)\mapsto(t,z-a(t,z))$,
we can assume that $\fff_\xi(\{u=u_0\})$ is horizontal.

\begin{figure}
    \centerline{\psfig{file=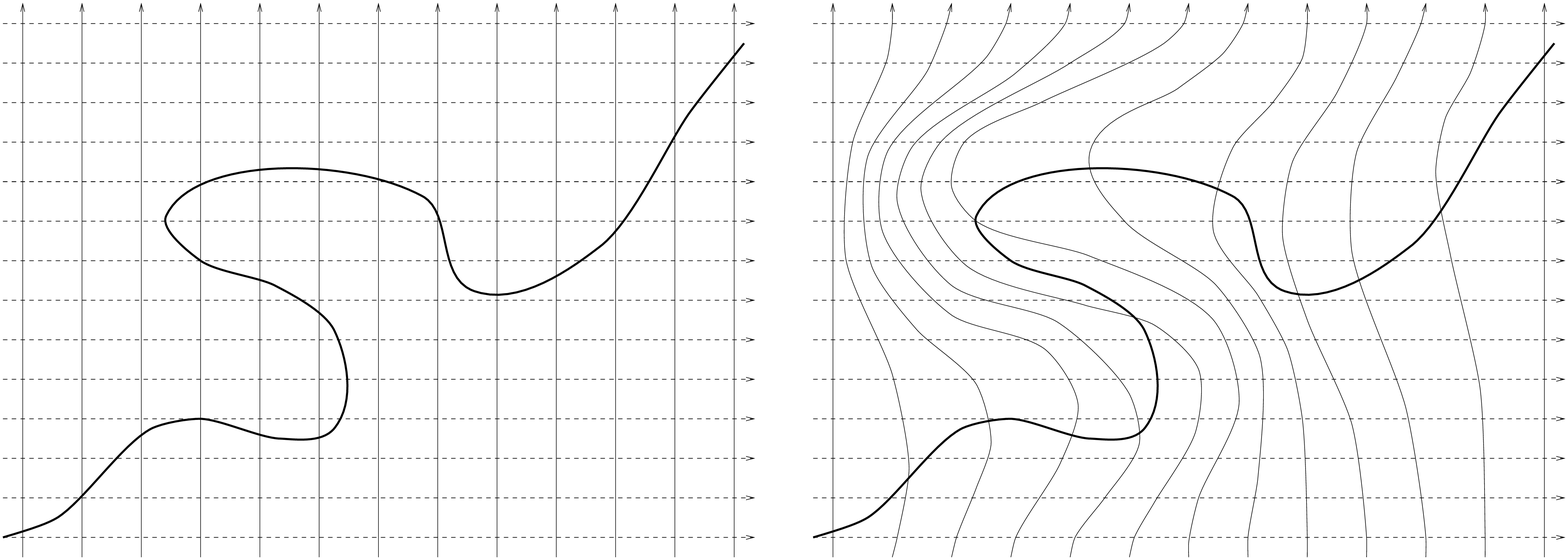,width=15.5truecm}}
    \mycaption{\label{newtrasv}
    Modification of the flow}
    \end{figure}

Now we know that the perturbed rectangle meets  
$\{u=u_0\}$ in a simple curve (not a graph any more, 
since we have changed coordinates).
Moreover, if one endows this curve and the leaves of 
$\fff_\xi(\{u=u_0\})$ with the correct orientation, one sees from
the proof of the elimination lemma and its inverse that 
negative tangencies never occur. Therefore the situation is
as in the left-hand side of Fig.~\ref{newtrasv},
and on the right-hand side of the same figure we suggest how to construct
the new flow $v_1$. This completes the proof.
\finedim{branched:elim:state}

We recall that, 
given a branched surface $P$, in~\cite{lnm} we have considered the 
(essentially unique) tangent vector field
to $P$ along $S(P)$ which always points to the right of $S(P)$
(i.e. it points from the locally 2-sheeted portion of $P$ to the locally
1-sheeted portion). Following the terminology of~\cite{christy}
we will call {\em maw} this field. 

\begin{prop}\label{normal:form:prop}
Let $P$ be a branched surface embedded in a contact $(N,\xi)$. Then:
\begin{enumerate}
\item Up to a $\cont0$-small perturbation of the embedding we can assume that  
$\fff_\xi(P)$ is S-stable  and directed by the maw at vertices;
\item If $N=\hatM$ and $P$ is a flow-spine which carries the homotopy class
of the field transversal to $\xi$ then up to isotopy (possibly not $\cont0$-small)
we can assume that $P$ is a faithful flow-spine for $\xi$ and $\fff_\xi(P)$ has
the same properties as in 1;
\item In both 1 and 2, if we start with an
$h$-embedding we can get an $h$-embedding.
\end{enumerate}
\end{prop}

\dim{normal:form:prop} To begin we note that in case 2 we can apply Proposition~\ref{homotopic:faithful}
(which involves a possibly non-small isotopy) and assume that $P$ is faithful for $\xi$. We will give now a
unified proof of all the statements, leaving to the reader the (easy) discussion on $h$-embeddings. 
This is because the transformations we will describe automatically preserve faithfulness.

\begin{figure}
    \centerline{\psfig{file=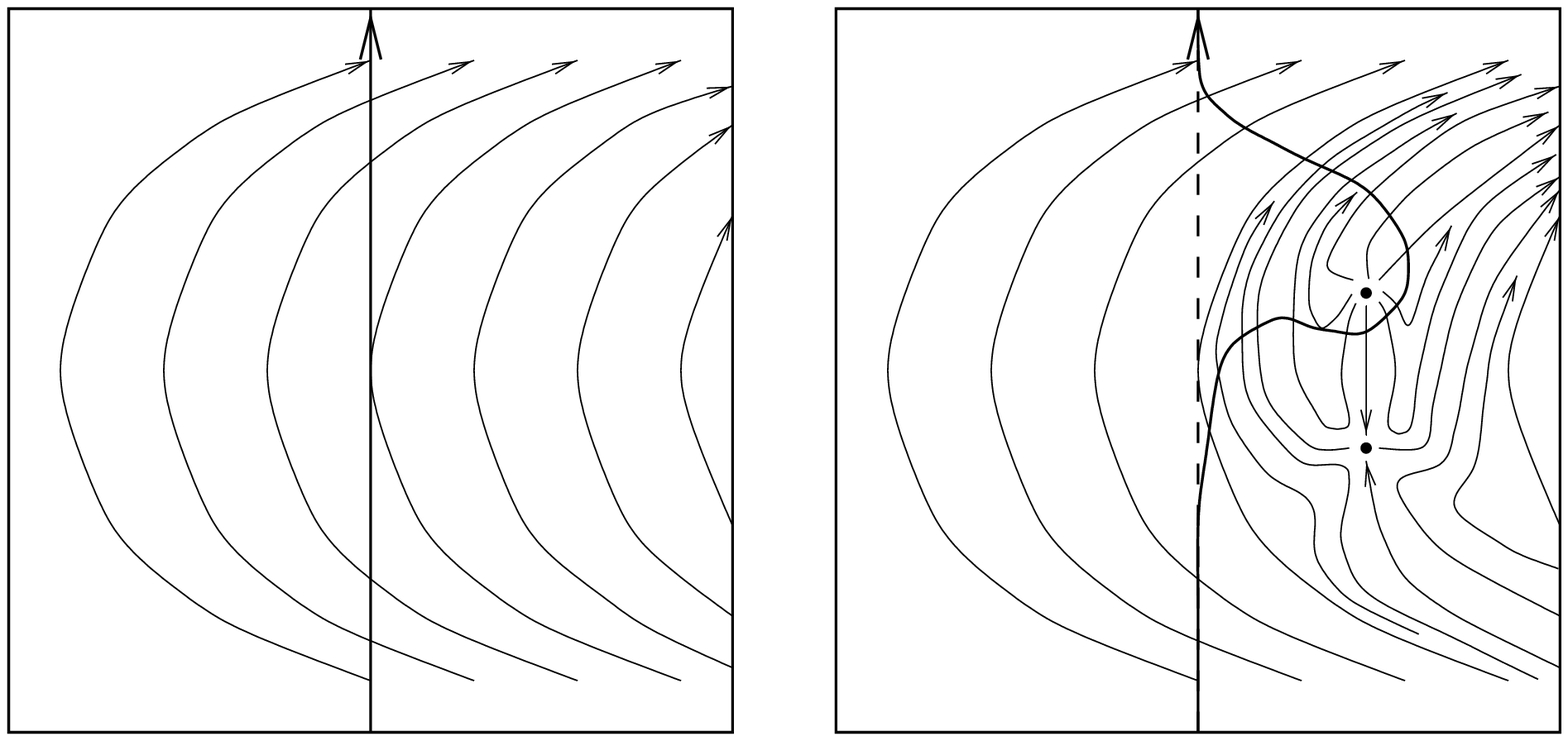,width=9truecm}}
    \mycaption{\label{posindex}
    How to remove a tangency of index $-1$}
    \end{figure}

First we modify the tangent plane to $P$ at vertices until
$\fff_\xi(P)$ is directed by  the maw there. Next, we take a small generic 
perturbation so that $\fff_\xi(P)$ has finitely many simple tangency points to $S(P)$. 
To conclude we remove tangencies of index $-1$ by first creating two positive singularities
(in a non-branched context) and then modifying the spine. A top view of what is happening
is shown in Fig.~\ref{posindex};
here one should imagine the portion which lies to the left of the thick line to
consist of two rectangles with foliations which are vertically aligned. The modification
of the spine is best understood as follows: imagine the single rectangle (on the right)
to be obtained by glueing an upper and a lower rectangle, and insert your finger from the left to
partially separate them. Note that Fig.~\ref{posindex} refers to one of the
two possible orientations for index $-1$, but the opposite orientation is dealt
with in a similar way.
\finedim{normal:form:prop}

\begin{prop}\label{all:close}
Let $P$ be a branched surface
with a $\cont\infty$-embedding $i_0$ in a contact manifold $(N,\xi)$.
Assume that $\ff_0=i_0^*(\fff_\xi(i_0(P)))$ is S-stable. 
Then any foliation $\ff_1$ on $P$ sufficiently $\cont\infty$ close to $\ff_0$
is induced by an embedding $i_1$ of $P$ in $N$ which is $\cont\infty$-close to $i_0$.
If $i_0$ is an $h$-embedding for some $h$ (or a faithful embedding 
of a flow-spine for $\xi$, or both) then the same holds for $i_1$. 
\end{prop}

\dim{all:close} 
Using stability near $S(P)$, as expressed in Proposition~\ref{S-stable:prop}, 
we can assume that $\ff_1$ coincides with $\ff_0$ near $S(P)$.

Let us identify for a moment $P$ with $i_0(P)$. Our first step
will be to find a contact structure $\xi'$ near $P$
which coincides with $\xi$ near $S(P)$, is arbitrarily
$\cont\infty$-close to $\xi$ and induces $\ff_1$ as characteristic
foliation on $P$. We only need to describe the extension of
$\xi$ across the discs of $P$, so the situation is as follows.
We have a disc $D$ embedded in a 3-manifold, a contact structure
$\xi$ near $D$, a foliation $\ff_1$ which is $\cont\infty$ close
to $\fff_\xi(D)$ on $D$ and coincides with $\fff_\xi(D)$ on a neighbourhood
$A$ (in $D$) of $\partial D$. We can trivialize a neighbourhood
of $D$ in $N$ as $D\times\mr$. Now $\xi$ can be viewed as a 1-form,
and $\fff_\xi(D)$ is defined by the pull-back $\alpha$ of $\xi$
with respect to the inclusion of $D$ in $N$. Therefore the foliation
$\ff_1$ is defined by a 1-form $\alpha+\omega$ on $D$, where
by assumption $\omega$ is $\cont\infty$ small and vanishes on $A$.
Using the coordinates $D\times\mr$ we can view $\omega$
as a form defined on the neighbourhood of $D$ (i.e.~$\omega$
is horizontal and vertically invariant). Now we can define
$\xi'$ as $\xi+\omega$. Of course $\xi'$ induces $\ff$ and
is $\cont\infty$ close to $\xi$, therefore it is a contact structure.
Moreover it coincides with $\xi$ near $\partial D$ because
$\omega$ vanishes on $A\times\mr$.

Having found $\xi'$, we define $\xi_s=(1-s)\xi+s\xi'$. 
Since $\xi'$ is close to $\xi$, this is a contact homotopy.
Moreover if we apply Moser's method the resulting time-depending 
vector field is constant near $S(P)$, and can be integrated up
to time 1 near $P$ yielding a diffeomorphism $f$ defined
on a neighbourhood of $P$, such that $f$ is $\cont\infty$ close
to the identity, $f$ is the identity near $S(P)$ and $f_*(\xi)=\xi'$.
Now it is sufficient to define $i_1=f^{-1}\compo i_0$ to conclude.
We leave to the reader the easy verification that if $i_0$ 
is an $h$-embedding or a faithful flow-spine then so is $i_1$.
\finedim{all:close}

\section{More uniqueness results}\label{more:uniqueness:sect}
In this section we prove the uniqueness parts of Theorems A and B
from the introduction. The reason for
not including them in Section~\ref{uniqueness:sect}
is that their proof uses the technique introduced to
establish Proposition~\ref{all:close}.

\begin{teo}\label{generic:determines}
Let $P$ be an abstract branched standard spine of $M$, 
let $\ff$ be an S-stable foliation on $P$
and let $i_j:P\to M$, $j=0,1$, be $\cont\infty$ embeddings. Then:
\begin{enumerate}
\item[(i)] If $\xi_0$ is a contact structure on $M$ which induces $(i_0)_*(\ff)$
on $i_0(P)$ then there exists a diffeomorphism
$\phi:M\to M$ such that $\phi_*(\xi_0)$ induces $(i_1)_*(\ff)$ on $i_1(P)$;
\item[(ii)] If $\xi_0$ and $\xi_1$ are contact structures on $M$ which induce
$(i_0)_*(\ff)$ and $(i_1)_*(\ff)$ on $i_0(P)$ and $i_1(P)$ respectively, then
there exist neighbourhoods $U_0$ and $U_1$ of $i_0(P)$ and $i_1(P)$
respectively and a contactomorphism $\phi:(U_0,\xi_0\ristr{U_0})\to (U_1,\xi_1\ristr{U_1})$;
\item[(iii)] If $\xi_0,\xi_1$ are as in (ii) and there exists a neighbourhood
of $i_0(P)$ on which $\xi_0$ is tight, then there exists a neighbourhood
of $i_1(P)$ on which $\xi_1$ is tight;
\item[(iv)] If $[\xi]$ denotes the isomorphism class of a 
contact structure $\xi$ on $\hatM=M\cup B$ then for 
$j=0,1$ the following sets coincide:
$\Big\{[\xi]:\ \fff_\xi(i_j(P))=(i_j)_*(\ff)\Big\}$.
\end{enumerate}
\end{teo}

It is perhaps useful, before the proof, to rephrase this result in less formal
terms: (i) means that the property for a foliation $\ff$ on $P$ of being induced
by a contact structure does not depend on the embedding of $P$; (ii) is
a uniqueness result which however, as already pointed out in the
introduction, cannot be interpreted in terms of germs;
(iii) shows that, even if the germ is not uniquely defined,
the property of it being tight is independent of the embedding of $P$;
(iv) means that the collection of contact structures
on $\hatM$ which induce $\ff$ on $P$ is again independent of the embedding.
Note that the diffeomorphism $\phi$ which appears in (i) and (ii), and tacitly
in (iv), does not map $i_0(P)$ to $i_1(P)$ in general.

\dim{generic:determines}
We start with (i). We first choose a diffeomorphism $f:M\to M$ such that
$f\compo i_1$ is arbitrarily $\cont\infty$-close to $i_0$ (for the existence 
of $f$ it is essential that $P$ be a standard spine). Therefore
$\ff=i_0^*(\fff_{\xi_0}(i_0(P)))$ and $\ff'=(f\compo i_1)^*(\fff_{\xi_0}((f\compo i_1)(P)))$
are arbitrarily $\cont\infty$-close together. Since $\ff$ is S-stable,
there exists a diffeomorphism $a:P\to P$ arbitrarily
$\cont\infty$-close to the identity such that $a_*(\ff')$ and $\ff$ coincide
on a neighbourhood of $S(P)$. Now one easily sees that the map
$$f\compo i_1\compo a\compo i_1^{-1}\compo f^{-1}: f(i_1(P))\to f(i_1(P))$$
extends to a diffeomorphism $h:M\to M$. By our choices we will have that
the foliations $(h\compo f\compo i_1)_*(\ff)$ and $\fff_{\xi_0}(f(i_1(P)))$
on $f(i_1(P))=h(f(i_1(P)))$ are $\cont\infty$-close and coincide
on a neighbourhood of the singular set. It follows, using the technique of
Proposition~\ref{all:close}, that we can modify the embedding $h\compo f\compo i_1$
away from $S(P)$ to an embedding $k:P\to M$ such that $\fff_{\xi_0}(k(P))=k_*(\ff)$.
Now, since the modification is $\cont\infty$-small and takes place away from the singular
set, it easily follows that there exists a diffeomorphism $\ell:M\to M$
such that $k=l\compo h\compo f \compo i_1$. The conclusion of (i) now follows
by taking $\phi=(\ell\compo h\compo f)^{-1}$.

To prove (ii) we only need to apply (i) and Theorem~\ref{foliation:determines}
to the branched surface $i_1(P)\subset M$ and the contact structures $\phi_*(\xi_0)$ and $\xi_1$.

Fact (iii) is now easy: if $(W,\xi_0\ristr{W})$ is tight then, up
to restricting $W$, we can assume that $W$ is diffeomorphic to $M$,
and apply (ii).

To prove (iv) we must show that given $\xi_0$ on $\hatM$ which induces 
$(i_0)_*(\ff)$ on $i_0(\ff)$ there exists an isomorphic $\xi_1$ which 
induces $(i_1)_*(\ff)$ on $i_1(\ff)$. To do this it is sufficient to extend
the map $\phi$ coming from (i) to a diffeomorphism of $\hatM$, which can be
done because $P$ is a standard spine, and define $\xi_1=\phi_*(\xi_0)$.
\finedim{generic:determines}

\begin{rem}
{\em To be completely formal in the above proof one should have given
a priori estimates on how close the foliations must be to be able
to apply the methods of Proposition~\ref{all:close} within $\xi_0$. 
Since one can restrict from the beginning to a compact neighbourhood of
$i_0(P)$ and prescribe a priori how close the various embeddings must be,
this is a technical point which can be safely left to the reader.}
\end{rem}

For Theorem B the following criterion is useful:

\begin{prop}\label{criterion}
Let $B$ be an open ball in a contact manifold $(\hatM,\xi)$.
Then $\xi$ is tight if and only if its restrictions to $B$ and to
a neighbourhood of $M\setminus B$ are tight.
\end{prop}

\dim{criterion}
The ``only if'' part is obvious. For the ``if'' part we first recall a general
definition and fact. Given $(\hatM,\xi)$ and $V\in\vecto(\hatM)$, we say that $V$ is
a contact field if the flow it generates leaves $\xi$ invariant. Of course this
definition makes sense also for partially defined vector fields. Now it is
a general fact~\cite{giroux:conv} that if $\alpha$ is a global equation of $\xi$
and $X$ is the Reeb field for $\alpha$ then for every $f\in\cont\infty(\hatM,\mr)$
there exists a unique $Y\in\vecto(\hatM)$ tangent to $\xi$ such that $f\cdot X+Y$
is a contact field. This implies quite easily that partially defined contact fields 
always extend to global ones.

Now assume $\xi$ is tight on $B$ and on its complement.  
Up to a small perturbation, we can assume that $\xi$ is tight near $\overline B$ and that
$\fff_\xi(\partial B)$ is the trivial foliation with one source, 
one sink, no saddles and no cycles. Therefore,
using~\cite{elia:palla}, we can identify $\overline B$ with the unit ball in the standard contact
structure $\xi_0=\d z-y\d x+x\d y$ on $\mr^3$. Now one easily checks that   
$x\cdot\partial/\partial x+y\cdot\partial/\partial y+2z\cdot\partial/\partial z$
is a contact field for $\xi_0$, therefore it extends to a contact field $V$ for $\xi$.
If by contradiction $\xi$ is overtwisted on $\hatM$ then there exists an
overtwisted disc $D$ which avoids $0\in B$. If $V$ generates $\{\phi_t\}$ then each $\phi_t$
is a contactomorphism of $\xi$, so $\fff_\xi(\phi_t(D))$ is ``constant''. Moreover
for $t$ big enough $\phi_t(D)$ is
contained in the complement of $B$, whence the contradiction.
\finedim{criterion}

\begin{cor}\label{tight:germ:def}
The property of carrying a tight contact structure on a  
neighbourhood is a well-defined property
of an abstract pair $(P,\ff)$ with S-stable $\ff$, 
and it is equivalent to carrying a global tight structure on $\hatM$. Moreover,
by Proposition~\ref{tight:determines}, 
such a structure on $\hatM$ is unique.
\end{cor}

\section{A smooth embedding theorem\\ and constructions of contact structures}\label{embed:cons:sect}

In this section we show that if $P$ is a branched standard
spine of a 3-manifold $M$ with boundary, then $M$ can be
embedded in $P\times\mr$ in a smooth fashion. We deduce from
this that the various contactization techniques known for
neighbourhoods of surfaces also work for branched spines, and
we apply these techniques to get contact structures which 
induce assigned characteristic foliations on a branched standard spine
(facts~\genericexists,~\positiveexistsclosed~and~\tightexistsflow~from 
the introduction). By simplicity we use the same notations as before,
but the results on the manifold with boundary $M$ hold with any boundary,
not necessarily $S^2$.

\begin{prop}\label{embed:prop}
Let $P$ be an embedded branched standard spine of a manifold $M$
with boundary. Then there exists a $\cont\infty$ embedding 
$i:M\to P\times\mr$ such that $i(P)$ is arbitrarily $\cont\infty$
close to $P\times\{0\}$.
\end{prop}

\dim{embed:prop}
As announced in the introduction, this result is a generalization of 
a theorem of Gillman and Rolfsen~\cite{g:r:1},~\cite{g:r:2} and a formal proof could be
given (with considerable effort) by modifying their explicit formulae
to get smooth functions. We prefer to describe the embedding pictorially: the
reader will be easily convinced that a formalization is indeed possible.

\begin{figure}
    \centerline{\psfig{file=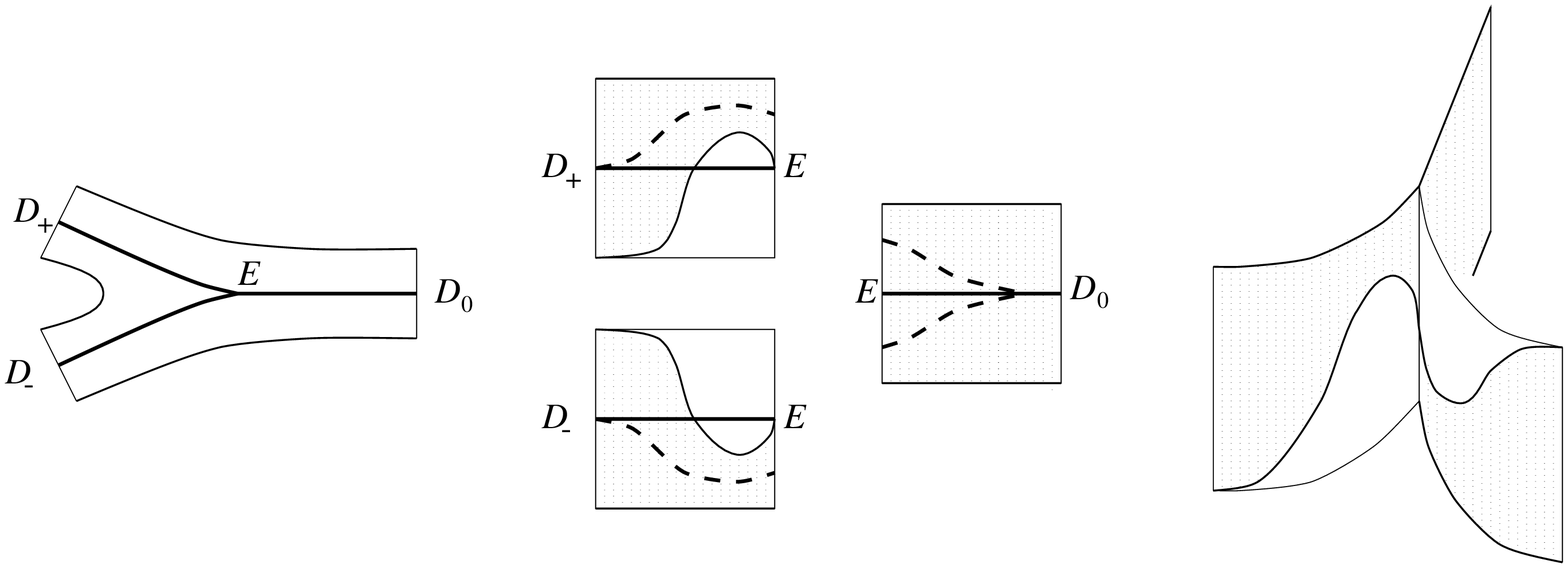,width=14truecm}}
    \mycaption{\label{embed1}
    Smooth embedding in two dimensions}
    \end{figure}

We start from a 2-dimensional situation, where a branched standard spine is a
train-track. In this case Fig.~\ref{embed1}
suggests how to proceed. On the left one sees a portion of
surface and the train-track embedded in it, and on the right the same
portion of surface is shown as a smooth subset of the product of the train-track 
with $\mr$. In the center we describe the same embedding avoiding 3-dimensional pictures,
and also showing the position of the train-track.

\begin{figure}
    \centerline{\psfig{file=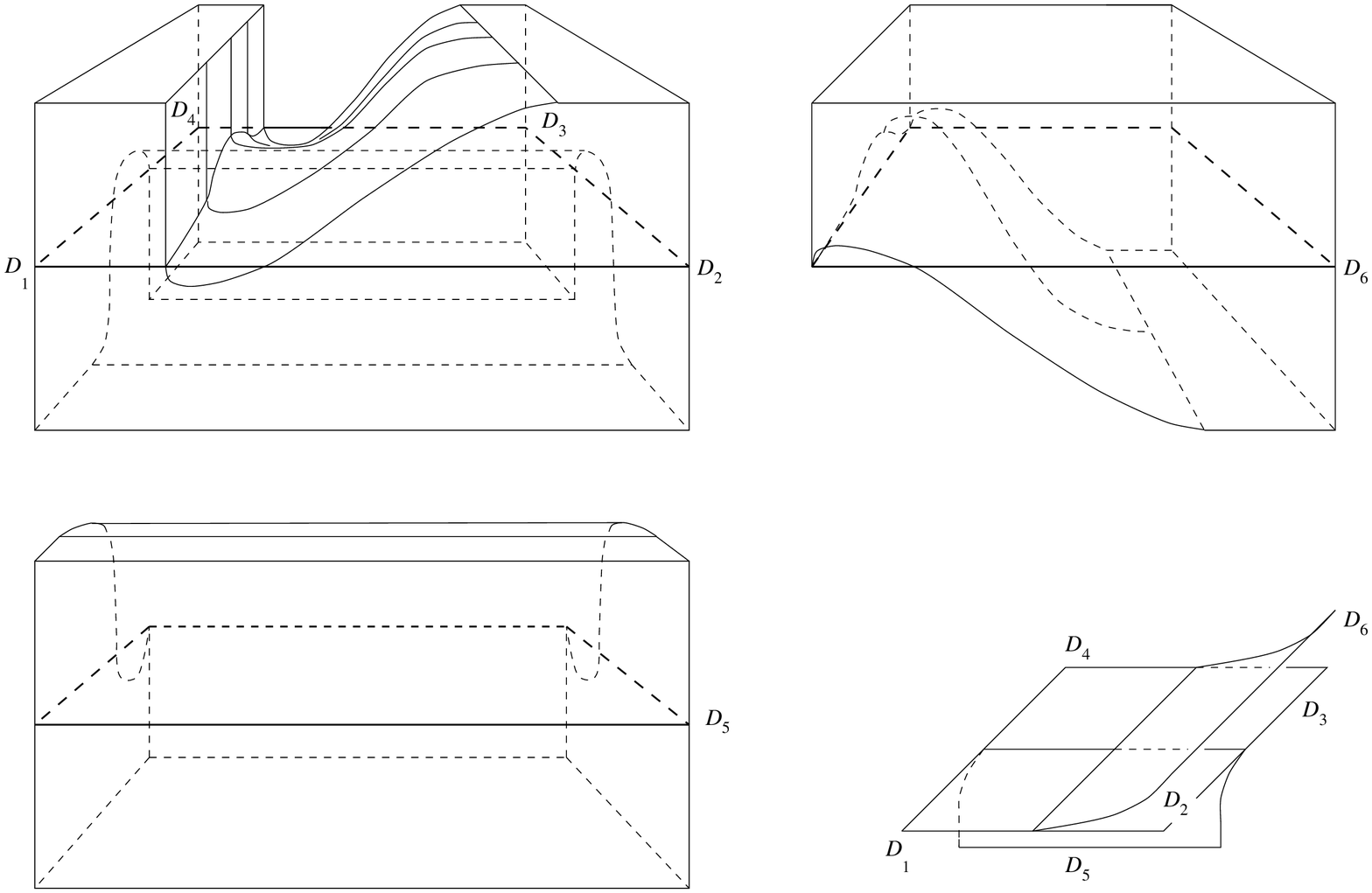,width=15.5truecm}}
    \mycaption{\label{embed2}
    Smooth embedding in three dimensions}
    \end{figure}

Now we go back to the 3-dimensional case. The subset of $P\times\mr$ diffeomorphic
to $M$ will contain $\{x\}\times[-1,1]$ for all $x\in P$ except near $S(P)$. 
Along singular edges, but
far from vertices, we only need to multiply the 2-dimensional picture by $\mr$.
The construction must be slightly modified near vertices, as shown in Fig.~\ref{embed2}.

The reader should note that in the 2-dimensional case the shaded regions of Fig.~\ref{embed1}
corresponding to $D_+$ and $D_-$ could be harmlessly interchanged at some
vertices. This cannot be
done in the 3-dimensional case, because one has to make
a coherent choice along edges and across vertices. 
Note also that to distinguish $D_+$ from $D_-$ we use
the screw-orientation of $P$ (i.e.~the fact that it is a spine of an oriented
manifold). 
\finedim{embed:prop}

\begin{prop}\label{non-zero-div:prop}
Let $P$ be a branched standard spine embedded in $M$, and let $\ff$ be an S-stable foliation
on $\ff(P)$ which has isolated singularities with 
non-zero divergence. Then there exists a contact structure $\xi$ on $M$ 
such that $\fff_\xi(P)=\ff$.
\end{prop}

\dim{non-zero-div:prop}
Our first step is to construct a contact structure $\xi$ on $P\times\mr$ which induces
the foliation $\ff$ on $P\times\{0\}$. Note that $P\times\mr$ is covered by
finitely many charts diffeomorphic to $\mr^3$ with smooth transition functions
(but the charts are not open), so it makes sense to speak of contact structures.
The definition of $\xi$ is exactly the same as in the case of surfaces~\cite{giroux:bourb}.
We take a form $\alpha$ which defines $\ff$ on $P$, we introduce on $P$
an area form $\omega$ and we consider the divergence $g$ of $\alpha$
with respect to $\omega$, characterized by the relation $\d\alpha=g\cdot\omega$.
Then we choose another form $\eta$ on $P$ with the same singularities as $\alpha$
and such that $\alpha\wedge\eta$ is a positive multiple of $\omega$
except at singularities.
Then, if $t$ is the coordinate on $\mr$, we just define $\xi$ as 
$\alpha+t(\d g-\eta)+g\d t$.

The second step (i.e.~the conclusion) consists once 
again in applying the techniques of
Proposition~\ref{all:close}. If we restrict $\xi$ from $P\times\mr$ to $M$
we have a contact structure which induces on a branched standard spine $P'$
arbitrarily $\cont\infty$ close to $P$ the desired characteristic foliation.
The same method used in the proof of Theorem~\ref{generic:determines}
allows to slightly isotope $\xi$  getting a new contact structure which induces
on $P$ exactly the desired foliation. Also in this case a complete formalization would 
require a priori estimates on closeness which we leave to the reader.
\finedim{non-zero-div:prop}

Now we are faced with the problem of extending the structure $\xi$ from $M$ to the
complementary ball $B=\hatM\setminus M$. Since $\xi$ extends as a plane field,
this can be deduced directly from the techniques of Eliashberg.
Alternatively, one can slightly enlarge
$B$ so that $\xi$ is defined on a neighbourhood of $\partial B$,  perturb 
$\partial B$ until $\fff_\xi(\partial B)$ becomes a Morse-Smale foliation and use
the following remark (inspired by~\cite{elia:over}) to be found in~\cite{giorgi}:

\begin{lem}
Let $\gg$ be a Morse-Smale foliation on $S^2$ and $\xi_0$ on $\mr^3$
be the standard overtwisted
structure. Then there exists an embedding $S^2\hookrightarrow\mr^3$ with
$\fff_{\xi_0}(S^2)=\gg$.
\end{lem}

This result, together with Proposition~\ref{non-zero-div:prop} proves 
fact~\genericexists~from the introduction, and therefore completes the proof of Theorem A.
We deal now with facts~\positiveexistsclosed~and~\tightexistsflow.
We start with a result concerning the manifold $M$ with boundary $S^2$,
later we will extend the structure to the ball $B$. Note that in the
next statement we use the obvious restriction from $\hatM$ to $M$
of the notion of faithful flow-spine.

\begin{prop}\label{positive-div:prop}
Let $P$ be a branched standard spine embedded in $M$, and let $\ff$ be an S-stable foliation
on $P$ which has isolated singularities with 
positive divergence. Then there exists a contact structure $\xi$ on $M$ 
such that $\fff_\xi(P)=\ff$ and $P$ is a faithful flow-spine for $\xi$.
\end{prop}

\dim{positive-div:prop}
Again we first define $\xi$ on $P\times\mr$. With the very same notations
as in the proof of Proposition~\ref{non-zero-div:prop} we define
$\xi=\alpha-t\cdot\eta+\d t$. Our choices and compactness of $P$ imply
quite easily that such a $\xi$ is indeed a contact structure on some
$P\times[-\varepsilon,\varepsilon]$. Now we can rescale the embedding
of $M$ in $P\times\mr$ so that it takes values in $P\times(-\varepsilon,\varepsilon)$.
Thus we have a contact structure $\xi$ on $M$, and $P$ is a faithful
flow-spine because the flow $\partial/\partial t$ in the coordinates
$P\times\mr$ is positively transversal to both $\xi$ and $P$, and the boundary
of $M$ has only one concave tangency curve.

Now the second step goes exactly as in Proposition~\ref{non-zero-div:prop}:
the property of $P$ being faithful for $\xi$ is
preserved under $\cont\infty$ small perturbations, whence the conclusion.
\finedim{positive-div:prop}

\begin{prop}\label{extend:homotopic}
Let $P,\ff,\xi$ be as in the previous proposition. Then $\xi$ extends to $\hatM$ 
within the homotopy class carried by $P$.
\end{prop}

\dim{extend:homotopic}
This fact can be proved in two ways. The most direct one is to use the
following result of Eliashberg: if $\xi$ is a plane field on $D^3$ and $\xi$
is contact near $S^2$ then there exists a contact structure $\xi'$ on $D^3$
homotopic to $\xi$ relatively to $S^2$. The second approach
consists in extending $\xi$ to an arbitrary contact structure and then
adjusting the homotopy class of $\xi$ using the methods of Lutz and Martinet
on the ball.
\finedim{extend:homotopic}

\begin{prop}\label{extend:tight:flow}
Let $P,\ff,\xi$ be as in Proposition~\ref{positive-div:prop}, and assume that 
$\xi$ is tight. Then $\xi$ extends to a tight structure on $\hatM$ for which
$P$ is a faithful flow-spine.
\end{prop}

\dim{extend:tight:flow}
The idea is to modify the boundary of the ball with the elimination lemma and then
fill the ball with the standard tight structure on the unit
ball in $\mr^3$. However one must be careful to preserve existence of a  
vector field which traverses the ball and is positively transversal to the structure.
If in the proof of Proposition~\ref{positive-div:prop} we start with a small
equation $\alpha$ of $\ff$ and in $\hatM$ we choose the sphere $\partial B$
very close to $P$, then we can imagine $\partial B$ as a very flat sphere,
e.g. one with equation $x^2+y^2+(z/r)^2=1$ for some very small $r>0$,
with traversing vector field close to $\partial/\partial z$ and 
contact structure close to $\d z$. This implies that the elimination lemma
is applied only within the upper and lower flat discs into which the
sphere splits, and it follows quite easily that the traversing field
can be adapted to the new position of the sphere. Now we can safely paste
the unit ball in the standard structure with its constant traversing vector field:
we only need to take a convex combination of the vector fields in a neighbourhood of the
boundary to make sure to get a global smooth vector field with the required
properties.
\finedim{extend:tight:flow} 

This completes the proof of Theorem B. 
 
\section{Concluding remarks}

According to Theorem B, if $P$ is a branched standard spine of  $\hatM$ and on $P$ we have an S-stable
foliation $\ff$  then the abstract pair $(P,\ff)$ carries at most one tight contact structure on $\hatM$
up to isomorphism.   The problem naturally arises to determine effectively which pairs $(P,\ff)$ carry
tight structures. The corresponding question for surfaces has been recently answered  by
Giroux~\cite{giroux:criteria}, under the (generically true) assumption that the foliation should admit a
splitting. Since the definition of splitting makes sense also for branched surfaces (where of course the
splitting curve will also be branched), one could speculate that the tightness of a sufficiently generic
$(P,\ff)$ depends on the topology of its splitting. However the question appears to be rather
challenging.

The difficulty in extending the tightness criterion to a branched context could have deep, not only
technical, reasons, possibly related to the following discussion. Let us first recall that for a foliated
surface $(\Sigma,\ff)$ with splitting curve $\Gamma$, Giroux's criterion is actually quite simple: the
germ of contact structure carried by $(\Sigma,\ff)$ is tight if and only either  $\Sigma=S^2$ and
$\Gamma$ is connected or $F\setminus\Gamma$ has no disc components. Since a neighbourhood of $\Sigma$
admits orientation-reversing automorphisms, this criterion is obviously the same for positive and for
negative contact structures. Forcing the analogy, assume that also the branched tightness criterion
(exists and) is formally the same for positive and for negative structures. Using the results of this
paper (in particular, existence) it would follow that a closed oriented 3-manifold $\hatM$ supports a
cooriented positive tight structure if and only if it supports a negative one. It was recently
established in~\cite{paolo} that there exist closed oriented 3-manifolds which carry simplectically
fillable (whence tight) positive contact structures, but do not carry any negative ones (for instance
the  Poincar\'e homology sphere).  As a consequence we would get examples of tight not symplectically
fillable structures.

Another very natural problem is to find necessary and/or sufficient conditions for two pairs $(P,\ff)$ and
$(P',\ff')$ to carry {\em the same} tight contact structure. In~\cite{lnm} we have provided a calculus
based on branched standard spines  for homotopy classes of plane fields, so one could hope to refine this
calculus to tight contact structures. Again, this does not seem to be straight-forward.  A solution of
these questions could probably be a significant contribution to the understanding of tight contact
structures on general 3-manifolds.

\vspace{.3cm}

\noindent\hspace{10.3cm}benedett@dm.unipi.it

\vspace{-3pt}

\noindent\hspace{10.3cm}petronio@dm.unipi.it

\vspace{3pt}

\noindent\hspace{10.3cm}Dipartimento di Matematica

\vspace{-3pt}

\noindent\hspace{10.3cm}Universit\`a di Pisa

\vspace{-3pt}

\noindent\hspace{10.3cm}Via F. Buonarroti, 2

\vspace{-3pt}

\noindent\hspace{10.3cm}I-56127 PISA, Italy

\end{document}